\documentclass[12pt]{article}
\font\smallit=cmti10

\usepackage{amsmath, amsthm, amssymb}
\begin{document} 

\begin{center}
{\bf ON PRIME RECIPROCALS IN THE CANTOR SET}
\vskip 20pt
{\bf Christian Salas}\\
{\smallit The Open University, Walton Hall, Milton Keynes MK7 6AA, United Kingdom}\\
{\tt c.p.h.salas@open.ac.uk}\\
\vskip 10pt
\end{center}
\vskip 30pt
\vskip 30pt 

\centerline{\bf Abstract}

\noindent
The middle-third Cantor set $\mathcal{C}_3$ is a fractal consisting of all the points in $[0, 1]$ which have non-terminating base-3 representations involving only the digits 0 and 2. It is easily shown that the reciprocals of all prime numbers $p > 3$ satisfying an equation of the form $2p + 1 = 3^q$ belong to $\mathcal{C}_3$. Such prime numbers have base-3 representations consisting of a contiguous sequence of 1's and are known as base-3 \emph{repunit} primes. It is natural to ask whether all prime numbers with reciprocals in $\mathcal{C}_3$ satisfy this equation. In this paper we show that the answer is no, but all primes with reciprocals in $\mathcal{C}_3$ do satisfy a closely related equation of the form $2pK + 1 = 3^q$. The base-3 repunit primes are thus shown to be a special case corresponding to $K = 1$. 

\pagestyle{myheadings}

\thispagestyle{empty} 
\baselineskip=15pt 
\vskip 30pt 

\section*{\normalsize 1 Introduction}

A prime number $p$ is called a base-$N$ repunit prime if it satisfies an equation of the form
\begin{equation}
\label{eq:DEF}
(N - 1)p + 1 = N^q
\end{equation} 
where $N \in \mathbb{N}-\{1\}$ and where $q$ is also prime. Such primes have the property that 
\begin{equation}
\label{eq:SUM}
p = \frac{N^q - 1}{N - 1} = \sum_{k=1}^q N^{q-k}
\end{equation}
so they can be expressed as a contiguous sequence of 1's in base $N$. For example, $p = 31$ satisfies \eqref{eq:DEF} for $N = 2$ and $q = 5$ and can be expressed as 11111 in base 2. The term \emph{repunit} was coined by A. H. Beiler \cite{BEIL} to indicate that numbers like these consist of repeated units. 

More importantly for what follows, the reciprocal of any such prime is an infinite series of the form
\begin{equation}
\label{eq:INFSUM}
\frac{1}{p} = \frac{N - 1}{N^q - 1} = \sum_{k=1}^{\infty} \frac{N-1}{N^{qk}}
\end{equation}
as can easily be verified using the usual methods for finding sums of series. Equation \eqref{eq:INFSUM} shows that $\frac{1}{p}$ can be expressed in base $N$ using only zeros and the digit $N - 1$. This single non-zero digit will appear periodically in the base-$N$ representation of $\frac{1}{p}$ at positions which are multiples of $q$. 
 
The case $N = 2$ corresponds to the famous Mersenne primes for which there are numerous important unsolved problems and a vast literature \cite{GUY}. They are sequence number A000668 in The Online Encyclopedia of Integer Sequences \cite{SLOANE}. The literature on base-$N$ repunit primes for $N \geq 3$ is principally concerned with computing and tabulating them for ever larger values of $N$ and $q$. An example is Dubner's \cite{DUB} tabulation for $2 \leq N \leq 99$ with large values of $q$. Relatively little is known about any peculiar mathematical properties that repunit primes in these other bases may possess.

In this paper we are prompted to investigate prime reciprocals belonging to the middle-third Cantor set $\mathcal{C}_3$ by the fact that $\mathcal{C}_3$ contains the reciprocals of all base-3 repunit primes, i.e., those primes $p$ which satisfy an equation of the form $2p + 1 = 3^q$ with $q$ prime. They are sequence number A076481 in The Online Encyclopedia of Integer Sequences. Simply putting $N = 3$ in \eqref{eq:DEF} and \eqref{eq:INFSUM} shows that $\frac{1}{p}$ can be expressed in base $3$ using only zeros and the digit $2$, which will appear periodically in the base-$3$ representation at positions which are multiples of $q$. Since only zeros and the digit $2$ appear in the ternary representation of $\frac{1}{p}$, $\frac{1}{p}$ is never removed in the construction of $\mathcal{C}_3$, so $\frac{1}{p}$ must belong to $\mathcal{C}_3$. In view of this, it is natural to then ask whether \emph{all} primes whose reciprocals belong to $\mathcal{C}_3$ satisfy an equation of the form $2p + 1 = 3^q$. In this paper we show that the answer to this question is no, but any prime number $p > 3$ whose reciprocal is in the middle-third Cantor set must satisfy a closely related equation of the form $2pK + 1 = 3^q$ where $q$ is not necessarily prime. The base-3 repunit primes are thus a special case corresponding to $K = 1$. 
\vskip 30pt
\section*{\normalsize 2 Preliminary definitions and results}

For easy reference in the discussion below, it is convenient to give a name to prime numbers whose reciprocals belong to $\mathcal{C}_3$. A logical one is the following:

{\bf Definition 1}
A Cantor prime is a prime number $p$ such that $\frac{1}{p} \in \mathcal{C}_3$.

We will also need to make use of some properties of prime numbers which result in the reciprocals of Cantor primes having a repeating cycle structure when represented in ternary form. The basic property is the following, which is commonly referred to as Fermat's Little Theorem:

{\bf Theorem 1} 
If $p$ is a prime and $a \neq p$ is another integer, then $a^{p-1} \equiv 1$ mod($p$). 

For a proof see, e.g., \cite{APOSTOL}. On the basis of Fermat's Little Theorem we can then define the order of an integer:

{\bf Definition 2}
If $p > 3$ is prime and $a \neq p$ is another integer, then the order of $a$ modulo $p$ is the least positive integer $q$ such that $a^q \equiv 1$ mod($p$).

The definition of the order of an integer assumes there is a unique least positive integer $q$ such that $a^q \equiv 1$ mod($p$). We can be sure of this on the basis of Fermat's Little Theorem since it guarantees that $q \leq p - 1$. 

In discussing Cantor primes we will be particularly concerned with the order of $3$ modulo $p$, i.e., the smallest integer $q$ such that $3^q \equiv 1$ mod($p$). Just as in the case of decimal representations of prime reciprocals, all ternary representations of prime reciprocals $\frac{1}{p}$ for $p > 3$ exhibit a repeating cycle which begins immediately after the point and has cycle length equal to the order of $3$ modulo $p$. To calculate the ternary representation of any prime reciprocal, say $\frac{1}{7}$, we can use a division algorithm as follows:
\begin{align*}
3 &= 0 \times 7 + 3 \\ 
9 &= 1 \times 7 + 2 \\ 
6 &= 0 \times 7 + 6 \\ 
18 &= 2 \times 7 + 4 \\ 
12 &= 1 \times 7 + 5 \\
15 &= 2 \times 7 + 1
\end{align*}
The remainder in each equation is multiplied by $3$ to obtain the left-hand side of the next equation and the procedure halts when a remainder of $1$ is obtained. From then on, the pattern of coefficients on the prime number will repeat indefinitely. In the above example a remainder of $1$ is obtained at the sixth step, so looking at the coefficients of $7$ we conclude that the ternary representation of $\frac{1}{7}$ is $0.\langle010212\rangle$ where the angle bracket notation is used to indicate that the string of digits in the bracket repeats indefinitely. This is completely analogous to the fact that the decimal, i.e., base 10, representation of $\frac{1}{7}$ has the repeating cycle structure $0.\langle142857\rangle$. The fact that the same thing happens in base 3 is crucial to obtaining the main result of this paper. Note that the cycle length is $6$ in the ternary representation of $\frac{1}{7}$ precisely because $6$ is the order of $3$ modulo $7$.
\vskip 30pt
\section*{\normalsize 3 Statement and proof of main result}

The following is a succinct statement of the main theorem we wish to prove in this paper:

{\bf Theorem 2} 
A prime number $p$ is a Cantor prime if and only if it satisfies an equation of the form $2pK + 1 = 3^q$ where $q$ is the order of 3 modulo $p$ and $K$ is a sum of non-negative powers of $3$ each of which is smaller than $3^q$. 

The base-3 repunit primes are then the special case in which $K = 3^0 = 1$. An example is 13, which satisfies $2p + 1 = 3^3$. In this special case a standard approach can be used to show that $q$ in $2p + 1 = 3^q$ must be prime if $p$ is prime \cite{OYE}. To see this, note that if $q = rs$ were composite we could obtain an algebraic factorisation of $3^q - 1$ as
\begin{equation}
3^q - 1 = (3^r)^s - (1)^s = (3^r - 1)(3^{(s-1)r} + 3^{(s-2)r} + \cdots + 1) 
\end{equation}  
We would then have
\begin{equation}
p = \frac{3^q - 1}{2} = \frac{(3^r - 1)}{2}(3^{(s-1)r} + 3^{(s-2)r} + \cdots + 1)
\end{equation}
Since $2\vert(3^r - 1)$, this would imply that $p$ is composite which is a contradiction. Therefore $q$ must be prime. In the more general case when $K \neq 1$, $q$ need not be prime. An example is 757, which satisfies $26p + 1 = 3^9$ with $K = 3^0 + 3^1 + 3^2 = 13$ and $q = 9$.

{\it Proof.}
In order to prove Theorem 2 it is necessary to consider the nature of $\mathcal{C}_3$ briefly. It is constructed recursively by first removing the open middle-third interval $(\frac{1}{3}, \frac{2}{3})$ from the closed unit interval $[0, 1]$. The remaining set is a union of two closed intervals $[0, \frac{1}{3}]$ and $[\frac{2}{3}, 1]$ from which we then remove the two open middle thirds $(\frac{1}{9}, \frac{2}{9})$ and $(\frac{7}{9}, \frac{8}{9})$. This leaves behind a set which is a union of four closed intervals from which we now remove the four open middle thirds, and so on. The set $\mathcal{C}_3$ consists of those points in $[0, 1]$ which are never removed when this process is continued indefinitely. 

Each $x \in \mathcal{C}_3$ can be expressed in ternary form as
\begin{equation}
\label{eq:TERN}
x = \sum_{k=1}^{\infty} \frac{a_k}{3^k} = 0.a_1a_2\ldots
\end{equation}
where all the $a_k$ are equal to 0 or 2. The construction of $\mathcal{C}_3$ amounts to systematically removing all the points in $[0, 1]$ which cannot be expressed in ternary form with only 0's and 2's, i.e., the removed points all have $a_k = 1$ for one or more $k \in \mathbb{N}$ \cite{OLM}. 

The construction of the Cantor set suggests some simple conditions which a prime number must satisfy in order to be a Cantor prime. If a prime number $p > 3$ is to be a Cantor prime, the first non-zero digit $a_{k_1}$ in the ternary expansion of $\frac{1}{p}$ must be 2. This means that for some $k_1 \in \mathbb{N}$, $p$ must satisfy 
\begin{equation} 
\label{eq:2P3PA}
\frac{2}{3^{k_1}} < \frac{1}{p} < \frac{1}{3^{k_1-1}}
\end{equation}
or equivalently
\begin{equation} 
\label{eq:2P3PB}
3^{k_1} \in (2p, 3p)
\end{equation}
Prime numbers for which there is no power of 3 in the interval $(2p, 3p)$, e.g., 5, 7, 17, 19, 23, 41, 43, 47, \ldots,\ can therefore be excluded immediately from further consideration. Note that there cannot be any other power of $3$ in the interval (2p, 3p) since $3^{k_1 - 1}$ and $3^{k_1 + 1}$ lie completely to the left and completely to the right of $(2p, 3p)$ respectively. 

If the next non-zero digit after $a_{k_1}$ is to be another 2 rather than a 1, it must be the case for some $k_2 \in \mathbb{N}$ that
\begin{equation}
\label{eq:P2PA}
\frac{2}{3^{k_1 + k_2}} < \frac{1}{p} - \frac{2}{3^{k_1}} < \frac{1}{3^{k_1 + k_2 - 1}} 
\end{equation}
or equivalently
\begin{equation}
\label{eq:P2PB}
3^{k_2} \in \bigg(\frac{2p}{3^{k_1} - 2p}, \frac{3p}{3^{k_1} - 2p}\bigg)
\end{equation}
Thus, any prime numbers which satisfy \eqref{eq:2P3PB} but for which there is no power of 3 in the interval $(\frac{2p}{3^{k_1} - 2p}, \frac{3p}{3^{k_1} - 2p})$ can again be excluded, e.g., 37, 113, 331, 337, 353, 991, 997, 1009. 

Continuing in this way, the condition for the third non-zero digit to be a 2 is
\begin{equation}
\label{eq:GENEQ0}
3^{k_3} \in \bigg(\frac{2p}{3^{k_2}(3^{k_1} - 2p) - 2p}, \frac{3p}{3^{k_2}(3^{k_1} - 2p) - 2p}\bigg)
\end{equation}
and the condition for the $n$th non-zero digit to be a 2 is
\begin{equation}
\label{eq:GENEQ1}
3^{k_n} \in \\ \\
\bigg(\frac{2p}{3^{k_{n-1}}(\cdots(3^{k_2}(3^{k_1} - 2p) - 2p)\cdots) - 2p}, \\ \\ \frac{3p}{3^{k_{n-1}}(\cdots(3^{k_2}(3^{k_1} - 2p) - 2p)\cdots) - 2p}\bigg)
\end{equation}
 
The ternary expansions under consideration are all non-terminating, so at first sight it seems as if an endless sequence of tests like these would have to be applied to ensure that $a_k \neq 1$ for any $k \in \mathbb{N}$. However, this is not the case: \eqref{eq:2P3PB} and \eqref{eq:GENEQ1} capture all the information that is required. To see this, let $p$ be a Cantor prime and let $3^{k_1}$ be the smallest power of 3 that exceeds $2p$. Since $p$ is a Cantor prime, both \eqref{eq:2P3PB} and \eqref{eq:GENEQ1} must be satisfied for all $n$. Multiplying \eqref{eq:GENEQ1} through by $3^{k_1-k_n}$ we get
\begin{equation}
\label{eq:GENEQ2}
3^{k_1} \in \bigg(\frac{3^{k_1-k_n}\cdot2p}{3^{k_{n-1}}(\cdots(3^{k_2}(3^{k_1} - 2p) - 2p)\cdots) - 2p}, \frac{3^{k_1-k_n}\cdot3p}{3^{k_{n-1}}(\cdots(3^{k_2}(3^{k_1} - 2p) - 2p)\cdots) - 2p}\bigg)
\end{equation}
Now, since all ternary representations of prime reciprocals $\frac{1}{p}$ for $p > 3$ have a repeating cycle which begins immediately after the point, it must be the case that $k_n = k_1$ for some $n$ in \eqref{eq:GENEQ2}. Setting $k_n = k_1$ in \eqref{eq:GENEQ2} we can therefore deduce from the fact that $3^{k_1} \in (2p, 3p)$ and the fact that \eqref{eq:GENEQ2} must be consistent with this for all values of $n$, that all Cantor primes must satisfy an equation of the form
\begin{equation}
\label{eq:GENEQ4}
3^{k_{n-1}}(\cdots(3^{k_2}(3^{k_1} - 2p) - 2p)\cdots) - 2p = 1
\end{equation}
where $k_1 + k_2 + \cdots + k_{n-1} = q$ is the cycle length in the ternary representation of $\frac{1}{p}$. By successively considering the cases in which there is only one non-zero term in the repeating cycle, two non-zero terms, three non-zero terms, etc., in \eqref{eq:GENEQ4}, and defining
\begin{align*}
d_1 &= q - k_1 \\
d_2 &= q - k_1 - k_2 \\
d_3 &= q - k_1 - k_2 - k_3 \\
\vdots \\
d_n &= q - k_1 - k_2 - \cdots - k_n = 0 
\end{align*}
it is easy to see that \eqref{eq:GENEQ4} can be rearranged as
\begin{equation}
\label{eq:GENEQ5}
2p\sum_{i=1}^n 3^{d_i} + 1 = 3^q
\end{equation} 
Setting $K = \sum_{i=1}^n 3^{d_i}$, we conclude that every Cantor prime must satisfy an equation of the form $2pK + 1 = 3^q$ as claimed in Theorem 2.

Conversely, every prime which satisfies an equation of this form must be a Cantor prime. To see this, note that we can rearrange \eqref{eq:GENEQ5} to get
\begin{equation}
\label{eq:GENEQ6}
\frac{1}{p} = \frac{2\sum_{i=1}^n 3^{d_i}}{3^q - 1} = 2\sum_{i=1}^n 3^{d_i}\bigg\{\frac{1}{3^q} + \frac{1}{3^{2q}} + \frac{1}{3^{3q}} + \cdots \bigg\}
\end{equation}  
Since $K = 2\sum_{i=1}^n 3^{d_i}$ involves only products of $2$ with powers of $3$ which are each less than $3^q$, \eqref{eq:GENEQ6} is an expression for $\frac{1}{p}$ which corresponds to a ternary representation involving only 2s. Thus, $\frac{1}{p}$ must be in the Cantor set if $2pK + 1 = 3^q$, and the proof is complete. 
\vskip 30pt 
\section*{\normalsize Acknowledgements}

I would like to thank anonymous referees who reviewed this manuscript.

\end{document}